\documentclass[reqno]{amsart}
\usepackage{amssymb,setspace}
\DeclareMathOperator{\arcsinh}{arcsinh}
\usepackage{ifpdf}
\ifpdf
 \usepackage[hyperindex,pagebackref]{hyperref}%
\else
 \expandafter\ifx\csname dvipdfm\endcsname\relax
 \usepackage[hypertex,hyperindex,pagebackref]{hyperref}%
 \else
 \usepackage[dvipdfm,hyperindex,pagebackref]{hyperref}%
 \fi
\fi
\theoremstyle{plain}
\newtheorem{thm}{Theorem}[section]
\numberwithin{equation}{section}
\allowdisplaybreaks[4]

\begin{document}

\title[sharp bounds for Neuman-S\'andor's mean]
{Sharp bounds for Neuman-S\'andor's mean in terms of the root-mean-square}

\author[W.-D. Jiang]{Wei-Dong Jiang}
\address[Jiang]{Department of Information Engineering, Weihai Vocational College, Weihai City, Shandong Province, 264210, China}
\email{\href{mailto: W.-D. Jiang <jackjwd@163.com>}{jackjwd@163.com}}

\author[F. Qi]{Feng Qi}
\address[Qi]{School of Mathematics and Informatics\\ Henan Polytechnic University\\ Jiaozuo City, Henan Province, 454010\\ China}
\email{\href{mailto: F. Qi <qifeng618@gmail.com>}{qifeng618@gmail.com}, \href{mailto: F. Qi <qifeng618@hotmail.com>}{qifeng618@hotmail.com}, \href{mailto: F. Qi <qifeng618@qq.com>}{qifeng618@qq.com}}
\urladdr{\url{http://qifeng618.wordpress.com}}



\begin{abstract}
In the paper, the authors find sharp bounds for Neuman-S\'andor's mean in terms of the root-mean-square.
\end{abstract}

\thanks{This work was partially supported by the Project of Shandong Province Higher Educational Science and Technology Program under grant No. J11LA57.}

\thanks{This paper was typeset using \AmS-\LaTeX}

\maketitle

\section{Introduction}

Throughout this paper, we assume that the real numbers $a$ and $b$ are positive and that $a\ne b$.
\par
The second Seiffert's mean $T(a,b)$ and Neuman-S\'andor's mean $M(a,b)$ are respectively defined in \cite{2, back10} by
\begin{equation}\label{seiffert-eq1.1}
 T(a,b)=\frac{a-b}{2\arctan \bigl(\frac{a-b}{a+b}\bigr)}\quad\text{and}\quad
 M(a,b)=\frac{a-b}{2\arcsinh \bigl(\frac{a-b}{a+b}\bigr)},
\end{equation}
while the arithmetic mean and the root-mean-square are respectively defined by
\begin{align}\label{seiffert-eq1.3}
A(a,b)=\frac{a+b}2\quad\text{and}\quad S(a,b)=\sqrt{\frac{a^2+b^2}2}\,.
\end{align}
A chain of inequalities between these four means
\begin{align*}
 A<M<T<S
\end{align*}
were established in~\cite{2,8}.
\par
In~\cite{back3}, the authors demonstrated that the double inequality
\begin{equation*}
\alpha S(a,b)+(1-\alpha)A(a,b)<T(a,b)<\beta S(a,b)+(1-\beta)A(a,b)
\end{equation*}
holds if and only if $\alpha\le \frac{4-\pi}{(\sqrt{2}\,-1)\pi}$ and $\beta\ge \frac{2}{3}$.
\par
In~\cite{4, 3}, the authors independently found that the double inequality
\begin{equation*}
\alpha S(a,b)+(1-\alpha)A(a,b)<M(a,b)<\beta S(a,b)+(1-\beta)A(a,b)
\end{equation*}
holds if and only if $\alpha\le \frac{1-\ln(1+\sqrt{2}\,)}{(\sqrt{2}\,-1)\ln(1+\sqrt{2}\,)}$ and $\beta\ge \frac{1}{3}$.
\par
In~\cite{seiffert10}, the authors discovered that the double inequality
\begin{equation*}
S\bigl(\alpha a+(1-\alpha)b,\alpha b+(1-\alpha)a\bigr)<T(a,b)
<S\bigl(\beta a+(1-\beta)b,\beta b+(1-\beta)a\bigr)
\end{equation*}
holds if and only if $\alpha\le\frac{1+\sqrt{16/\pi^2-1}\,}2$ and $\beta\ge\frac{3+\sqrt6\,}6$.
\par
Motivated by the above double inequalities, we naturally ask a question: What are the best constants $\alpha\ge\frac12$ and $\beta\le1$ such that the double inequality
\begin{equation}\label{seiffert-eq1.4}
S\bigl(\alpha a+(1-\alpha)b,\alpha b+(1-\alpha)a\bigr)<M(a,b)
<S\bigl(\beta a+(1-\beta)b,\beta b+(1-\beta)a\bigr)
\end{equation}
holds for all $a,b>0$ with $a\ne b$?
\par
The aim of this paper is just to give an affirmative answer to this question.
\par
The main result of this paper may be formulated as the following theorem.

\begin{thm}\label{seiffert-th1.1}
The double inequality~\eqref{seiffert-eq1.4} holds true if and only if
\begin{equation*}
\alpha\le \frac12\Biggl\{1+\sqrt{\frac{1}{\bigl[\ln \bigl(1+\sqrt{2}\,\bigr)\bigr]^2}-1}\,\Biggr\}=0.76\dotsc\quad \text{and}\quad \beta\ge \frac{3+\sqrt{3}\,}{6}=0.78\dotsc.
\end{equation*}
\end{thm}

\section{Proof of Theorem~\ref{seiffert-th1.1}}

For simplicity, denote
\begin{equation*}
\lambda=\frac12\Biggl\{1+\sqrt{\frac{1}{\bigl[\ln \bigl(1+\sqrt{2}\,\bigr)\bigr]^2}-1}\,\Biggr\}\quad \text{and}\quad \mu=\frac{3+\sqrt{3}\,}{6}.
\end{equation*}
It is clear that, in order to prove the double inequality~\eqref{seiffert-eq1.4}, it suffices to show
\begin{equation}\label{seiffert-eq2.1}
 M(a,b)>S\bigl(\lambda a+(1-\lambda)b,\lambda b+(1-\lambda)a\bigr)
\end{equation}
and
\begin{equation}\label{seiffert-eq2.2}
 M(a,b)<S\bigl(\mu a+(1-\mu)b,\mu b+(1-\mu)a\bigr).
\end{equation}
\par
From definitions in~\eqref{seiffert-eq1.1} and~\eqref{seiffert-eq1.3}, we see that both $M(a,b)$ and $S(a,b)$ are symmetric and homogeneous of degree $1$. Hence, without loss of generality, we assume that $a>b>0$. If replacing $\frac{a}b>1$ by $t>1$ and letting $p\in \bigl(\frac12,1\bigr)$, then
\begin{multline}\label{seiffert-eq2.3}
S(pa+(1-p)b,pb+(1-p)a)-M(a,b)\\
=\frac{b}2\frac{\sqrt{[pt+(1-p)]^2+[p+(1-p)t]^2}\,}{\arcsinh\frac{t-1}{t+1}}f(t),
\end{multline}
where
\begin{equation}\label{seiffert-eq2.4}
f(t)=\sqrt{2}\,\arcsinh \frac{t-1}{t+1}-\frac{t-1}{\sqrt{[pt+(1-p)]^2+[p+(1-p)t]^2}\,}.
\end{equation}
Standard computations lead to
\begin{align}
f(1)&=0,\label{seiffert-eq2.5}\\
\label{seiffert-eq2.6}
 \lim_{t\to\infty}f(t)&=\sqrt{2}\,\ln\bigl(1+\sqrt{2}\,\bigr)-\frac{1}{\sqrt{2p^2-2p+1}\,},
\end{align}
and
\begin{equation}\label{seiffert-eq2.7}
 f'(t)=\frac{f_1(t)}{(1+t)\sqrt{1+t^2}\,\{[pt+(1-p)]^2+[p+(1-p)t]^2\}^{3/2}},
\end{equation}
where
\begin{equation}\label{seiffert-eq2.8}
 f_1(t)=2\bigl\{[pt+(1-p)]^2+[p+(1-p)t]^2\bigr\}^{3/2}-(1+t)^2\sqrt{1+t^2}\,
\end{equation}
and
\begin{equation}\label{seiffert-eq2.9}
\Bigl\{2\bigl([pt+(1-p)]^2+[p+(1-p)t]^2\bigr)^{3/2}\Bigr\}^2 -\Bigr[(1+t)^2\sqrt{1+t^2}\,\Bigl]^2=(t-1)^2g_1(t)
\end{equation}
with
\begin{equation}
\begin{split}\label{seiffert-eq2.10}
g_1(t)&=\bigl(32p^6-96p^5+144p^4-128p^3+72p^2-24p+3\bigr)t^4\\
&\quad-2\bigl(64p^6-192p^5+240p^4-160p^3+48p^2-1\bigr)t^3\\
&\quad+6\bigl(32p^6-96p^5+112p^4-64p^3+24p^2-8p+1\bigr)t^2\\
&\quad-2\bigl(64p^6-192p^5+240p^4-160p^3+48p^2-1\bigr)t\\
&\quad+32p^6-96p^5+144p^4-128p^3+72p^2-24p+3
\end{split}
\end{equation}
and
\begin{align}\label{seiffert-eq2.11}
 g_1(1)=16\bigl(6p^2-6p+1\bigr).
\end{align}
\par
Let
\begin{equation*}
g_2(t)=\frac{g_1'(t)}2,\quad g_3(t)=\frac{g_2'(t)}6,\quad\text{and}\quad g_4(t)=g_3'(t).
\end{equation*}
Then simple computations result in
\begin{align}
\begin{split}\label{seiffert-eq2.12}
g_2(t)&=2\bigl(32p^6-96p^5+144p^4-128p^3+72p^2-24p+3\bigr)t^3\\
&\quad-3\bigl(64p^6-192p^5+240p^4-160p^3+48p^2-1\bigr)t^2\\
&\quad+6\bigl(32p^6-96p^5+112p^4-64p^3+24p^2-8p+1\bigr)t\\
&\quad-\bigl(64p^6-192p^5+240p^4-160p^3+48p^2-1\bigr),
\end{split}\\
\label{seiffert-eq2.13}
 g_2(1)&=16\bigl(6p^2-6p+1\bigr),\\
\begin{split}\label{seiffert-eq2.14}
 g_3(t)&=\bigl(32p^6-96p^5+144p^4-128p^3+72p^2-24p+3\bigr)t^2\\
 &\quad-\bigl(64p^6-192p^5+240p^4-160p^3+48p^2-1\bigr)t\\
 &\quad+32p^6-96p^5+112p^4-64p^3+24p^2-8p+1,
\end{split}\\
\label{seiffert-eq2.15}
 g_3(1)&=16p^4-32p^3+48p^2-32p+5,\\
\begin{split}\label{seiffert-eq2.16}
 g_4(t)&=2\bigl(32p^6-96p^5+144p^4-128p^3+72p^2-24p+3\bigr)t\\
 &\quad-\bigl(64p^6-192p^5+240p^4-160p^3+48p^2-1\bigr),
\end{split}\\
\label{seiffert-eq2.17}
 g_4(1)&=48p^4-96p^3+96p^2-48p+7.
\end{align}
\par
When $p=\lambda$, the quantities~\eqref{seiffert-eq2.6}, \eqref{seiffert-eq2.11}, \eqref{seiffert-eq2.13}, \eqref{seiffert-eq2.15}, and~\eqref{seiffert-eq2.17} become
\begin{gather}\label{seiffert-eq2.18}
 \lim_{t\to\infty} f(t)=0,\\\label{seiffert-eq2.19}
 g_1(1)=g_2(1)=-\frac{8\bigl\{4\bigl[\ln\bigl(1+\sqrt{2}\,\bigr)\bigr]^2-3\bigr\}} {\bigl[\ln\bigl(1+\sqrt{2}\,\bigr)\bigr]^2}<0,\\
\label{seiffert-eq2.21}
 g_3(1)=-\frac{7\bigl[\ln\bigl(1+\sqrt{2}\,\bigr)\bigr]^4 -4\bigl[\ln\bigl(1+\sqrt{2}\,\bigr)\bigr]^2-1} {\bigl[\ln\bigl(1+\sqrt{2}\,\bigr)\bigr]^4}<0,\\\label{seiffert-eq2.22}
 g_4(1)=-\frac{5\bigl[\ln\bigl(1+\sqrt{2}\,\bigr)\bigr]^4-3} {\bigl[\ln\bigl(1+\sqrt{2}\,\bigr)\bigr]^4}<0,
\end{gather}
and
\begin{equation}\label{seiffert-eq2.23}
32p^6-96p^5+144p^4-128p^3+72p^2-24p+3
=-\frac{2\bigl[\ln\bigl(1+\sqrt{2}\,\bigr)\bigr]^6-1}{\bigl[2\ln\bigl(1+\sqrt{2}\,\bigr)\bigr]^6}>0.
\end{equation}
Consequently, from~\eqref{seiffert-eq2.10}, \eqref{seiffert-eq2.12}, \eqref{seiffert-eq2.14}, \eqref{seiffert-eq2.16}, and~\eqref{seiffert-eq2.18}, it is very easy to obtain that
\begin{equation}\label{seiffert-eq2.24}
 \lim_{t\to\infty} g_1(t)=\infty,\quad
 \lim_{t\to\infty} g_2(t)=\infty,\quad
 \lim_{t\to\infty} g_3(t)=\infty,\quad
 \lim_{t\to\infty} g_4(t)=\infty.
\end{equation}
From~\eqref{seiffert-eq2.16} and~\eqref{seiffert-eq2.23}, it is immediate to derive that the function $g_4(t)$ is strictly increasing on $[1,\infty)$, and so, by virtue of~\eqref{seiffert-eq2.22} and the final limit in~\eqref{seiffert-eq2.24}, there exists a point $t_0>1$ such that $g_4(t)<0$ on $[1,t_0)$ and $g_4(t)>0$ on $(t_0,\infty)$. Hence, the function $g_3(t) $ is strictly decreasing on $[1,t_0]$ and strictly increasing on $[t_0,\infty)$.
Similarly, by~\eqref{seiffert-eq2.21} and the third limit in~\eqref{seiffert-eq2.24}, there exists a point $t_1>t_0>1$ such that $g_2(t)$ is strictly decreasing on $[1,t_1]$ and strictly increasing on $[t_1,\infty)$. Further, by~\eqref{seiffert-eq2.19} and the second limit in~\eqref{seiffert-eq2.24}, there exists a point $t_2>t_1>1$ such that $g_1(t)$ is strictly decreasing on $[1,t_2]$ and strictly increasing on $[t_2,\infty)$.
Thereafter, by~\eqref{seiffert-eq2.7} to~\eqref{seiffert-eq2.9}, \eqref{seiffert-eq2.19}, and the first limit in~\eqref{seiffert-eq2.24}, there exists a point $t_3>t_2>1$ such that $f(t)$ is strictly decreasing on $[1,t_3]$ and strictly increasing on $[t_3,\infty)$.
As a result, the inequality~\eqref{seiffert-eq2.1} follows from equations~\eqref{seiffert-eq2.3} to~\eqref{seiffert-eq2.5}, and~\eqref{seiffert-eq2.18}, together with the piecewise monotonicity of $f(t)$.
\par
When $p=\mu$, the equation~\eqref{seiffert-eq2.10} becomes
\begin{equation}\label{seiffert-eq2.28}
 g_1(t)=\frac{5t^2+8t+5}{27}(t-1)^2>0
\end{equation}
for $t>1$. By equations~\eqref{seiffert-eq2.7} to~\eqref{seiffert-eq2.10} and the inequality~\eqref{seiffert-eq2.28}, it can be concluded that $f(t)$ is strictly increasing and positive on $[1,\infty)$. The inequality~\eqref{seiffert-eq2.2} follows.
\par
It is not difficult to verify that the mean $S\bigl(xa+(1-x)b,xb+(1-x)a\bigr)$ is continuous and strictly increasing on $\bigl[\frac12,1\bigr]$.
From this monotonicity and inequalities~\eqref{seiffert-eq2.1} and~\eqref{seiffert-eq2.2}, one can conclude that the double inequality~\eqref{seiffert-eq1.4} holds true for all $\alpha\le\lambda$ and $\beta\ge \mu$.
\par
For any given number $p$ satisfying $1>p>\lambda$, it is obvious that the limit~\eqref{seiffert-eq2.6} is positive. This positivity, together with~\eqref{seiffert-eq2.3} and~\eqref{seiffert-eq2.4}, implies that for $1>p>\lambda$ there exists $T_0=T_0(p)>1$ such that the inequality
\begin{equation*}
 S(pa+(1-p)b,pb+(1-p)a)>M(a,b)
\end{equation*}
holds for $\frac{a}b\in (T_0,\infty)$. This tells us that the constant $\lambda$ is the best possible.
\par
For $\frac12<p<\mu$, from~\eqref{seiffert-eq2.11} one has
\begin{align}\label{seiffert-eq2.29}
 g_1(1)=16(6p^2-6p+1)<0.
\end{align}
From the inequality~\eqref{seiffert-eq2.29} and the continuity of $g_1(t)$, there exists a number $\delta=\delta(p)>0$ such that the function $g_1(t)$ is negative on $(1,1+\delta)$. This negativity, together with~\eqref{seiffert-eq2.3}, \eqref{seiffert-eq2.5}, \eqref{seiffert-eq2.7}, and~\eqref{seiffert-eq2.10}, implies that for any $\frac12<p<\mu$, there exists $\delta=\delta(p)>0$ such that the inequality
\begin{equation*}
 M(a,b)>S(pa+(1-p)b,pb+(1-p)a)
\end{equation*}
is valid for $\frac{a}b\in (1,1+\delta)$. Consequently, the number $\mu$ is the best possible. The proof of Theorem~\ref{seiffert-th1.1} is complete.

\end{document}